\documentclass[12pt]{article}
\usepackage[mathscr]{eucal}
\usepackage{amsbsy}
\usepackage{amssymb}
\usepackage{amsmath}
\usepackage{amsthm}
\usepackage{graphics}
\usepackage{epsfig}

\usepackage{latexsym}

\newcommand{\PP}{\mathbf{P}}

\def\phi{\varphi}
\def\bbn{{\mathbb N}}
\def\bbr{{\mathbb R}}

\def\bz{{\mathbb Z}} 
\def\ba{{\bf a}}
\def\bb{{\bf b}}
\def\bc{{\bf c}}
\def\bg{{\bf g}}
\def\bh{{\bf h}}
\def\be{{\bf e}}
\def\bg{{\bf g}}
\def\bbf{{\bf f}}
\def\bbh{{\bf H}}
\def\bbs{{\bf S}}
\def\bll{{\bf l}}

\def\bk{{\bf k}}
\def\bt{{\bf t}}
\def\bs{{\bf s}}
\def\bu{{\bf u}}
\def\bv{{\bf v}}
\def\bll{{\bf l}}
\def\bx{{\bf x}}
\def\bbx{{\bf X}}
\def\bY{{\bf Y}}
\def\bT{{\bf T}}
\def\bn{{\bf n}}

\def\bi{{\bf i}}
\def\Bl1{{\bf 1}}
\def\B2{{\bf 2}}
\def\B0{{\bf 0}}
\def\by{{\bf y}}

\def\a{\alpha}
\def\b{\beta}

\def\vfi{\varphi}
\def\e{\varepsilon}
\def\g{\gamma}
\def\G{\Gamma}
\def\l{\lambda}

\def\cl{{\mathscr L}}

\def\H{{\mathscr H}}
\def\C{{\mathscr C}}
\def\U{{\mathscr U}}
\def\=A8{\"o}

\def\lygfdd{\stackrel{\rm d}{=}}

\def\fdd{\stackrel{f.d.d.}{\longrightarrow}}

\newcommand{\beq}{\begin{equation}}
\newcommand{\eeq}{\end{equation}}
\newcommand\beqn{\begin{displaymath}}  
\newcommand\eeqn{\end{displaymath}}

\newcommand{\halmos}{\vspace{3mm} \hfill \mbox{$\Box$}\\[2mm]}
\theoremstyle{plain}
\newtheorem{teo}{Theorem}

\newtheorem{cor}[teo]{Corollary}
\newtheorem{prop}[teo]{Proposition}
\theoremstyle{definition}
\newtheorem{definition}[teo]{Definition}

\begin{document}

\title{ Lamperti type theorems for random fields
\footnotemark[0]\footnotetext[0]{ \textit{Short title:} Lamperti theorems }
\footnotemark[0]\footnotetext[0]{%
\textit{MSC 2000 subject classifications}. Primary 60G18, secondary
60G60, 26B99.} \footnotemark[0] \footnotemark[0]\footnotetext[0]{ \textit{Corresponding
author:} Vygantas Paulauskas, Vilnius University, Department
of Mathematics and Informatics,
 e-mail: vygantas.paulauskas@mif.vu.lt  }}

\author{ Youri Davydov$^{\text{\small 1}}$  and  Vygantas Paulauskas$^{\text{\small 2}}$ \\
{\small $^{\text{1}}$ Universit\'e Lille 1, Laboratoire Paul
Painlev\'e and St. Petersburg University}\\
 {\small $^{\text{2}}$ Vilnius University, Department
of Mathematics and Informatics}
 }

\maketitle

\begin{abstract}
In the paper we consider Lamperti type theorems for random fields. Together with known results we present some new results
on $\bbr^m$-valued self-similar fields $\{\bbx (\bt), \ \bt \in \bbr^d \}$, their domains of attraction and the so-called Lamperti transformations, expressing the relation between self-similarity
and stationarity. Also we investigate regularly and slowly varying functions of several variables.
\end{abstract}
\section{Introduction }

It is well-known what important role in the theory of stochastic processes plays the so-called self-similar processes, which were introduced by J. Lamperti  in his fundamental paper \cite{Lamperti} (only he used the term semi-stable, the name self-similar was introduced by B.B. Mandelbroit later). Since the formulations of the main results from \cite{Lamperti} are quite simple, we provide them here, only we use the word self-similar instead of semi-stable and letter $H$ (instead of $\a$) for the index of self-similarity.
For two  processes $\bbx(t)$ and $\bY(t)$, \ $t\in T$, with values in $\bbr^m$, we write $\{\bbx(t),\ t\in T\}\stackrel{\rm d}{=} \{\bY(t), \ t\in T\}$ if the finite-dimensional distributions (f.d.d.) coincide.  A random vector is called full if it is not supported on a lower dimensional hyperplane. If a set $T$ is $\bbr^d$ or some subset of $\bbr^d$ then a process $\{\bbx(t), \bt\in T\}$ is called proper, if $\bbx(t)$ it full for every $\bt \ne 0$. Here and in what follows letters in bald stand for vectors, and equalities or inequalities between them are componentwise. If $\bt, \bs \in \bbr^d$, then we denote $\bt\cdot \bs=(b_1s_1, \dots, b_ds_d)$ and, if $t_i\ne 0, 1\le i\le d,$ \ $\bt^{-1}=(t_1^{-1}, \dots , t_d^{-1}).$ A diagonal matrix with entries $a_i$ on the diagonal will be denoted by ${\rm diag}(a_1, \dots , a_d)$.

\begin{definition}\label{deflamperti} (\cite{Lamperti}) A $\bbr^m$-valued process $\{\bbx(t), t\in \bbr_+\}$ is called self-similar, if it is proper and stochastically continuous and for any $a>0$ there exist constants $b(a)$ and $\bc(a)$ such that
$$\{\bbx(at), \ t\in \bbr_+\}\stackrel{\rm d}{=} \{b(a)\bbx(t)+\bc(a), \ t\in \bbr_+\}$$.
\end{definition}

If we compare this definition with those, which were used later (see, for example Def. 7.1.1 in \cite{Samorod}), we see that the  requirement of stochastic continuity and the shifts are dropped (i.e., $\bc(a)\equiv 0$). The main results of \cite{Lamperti} are the following two theorems.
\begin{teo}\label{thm1} (\cite{Lamperti}) If $\{ \bbx(t),, \ t\in \bbr_+\} $ is self-similar, then there exists $H \ge 0$ such that
$$
b(a)=a^H.
$$
If $H>0$, then the distribution of $X(0)$ is concentrated at a point $\omega\in \bbr^m$ and
$$
\bc(a)=\omega  (1-a^H).
$$
If $H=0$, then $\bc(a)\equiv 0$ and the process is trivial in the sense that $\bbx(t)=\bbx(0) \ a.s.$ for each $t$.
\end{teo}
In (\cite{Lamperti}) the parameter $H$ was called the order of the self-similar process, nowadays often it is called Hurst index of self-similarity.

\begin{teo}\label{thm2} (\cite{Lamperti}) Let $\{ \bY(t)\}$ be a $\bbr^m$-valued stochastic process such that there exist a real valued function $f(a)$, $0< f(a) \to \infty,$ as $a \to \infty,$ $\bbr^m$-valued function $\bk(a)$, and a proper stochastically continuous process $\bbx(t)$, satisfying the relation
\begin{equation}\label{lim1}
\lim_{a \to \infty}\left \{ \frac{\bY(a t)+\bk(a)}{f(a)}\right \} \fdd \{ \bbx(t)\},
\end{equation}
where $\fdd$ stands for the  convergence of f.d.d.. Then the process $\{\bbx(t)\}$ is self-similar and
$$
f(a)=a ^H L(a), \quad \bk(a)=\omega(a)a ^H L(a),
$$
where $H>0$, $L(a)$ is a slowly varying function and $\lim_{a \to \infty}\omega (a)=w, \ w\in \bbr^m$. The order of $\{\bbx(t)\}$ is $H$ and $\bbx(0)=w.$ Conversely, every self-similar process $\{\bbx(t)\}$  of positive order can be obtained as in (\ref{lim1}) for some process $\{\bY(t)\}$.
\end{teo}

One more result, which was not considered important by Lamperti himself, since it was not separated as some statement, but which later was named as Lamperti transformation, can be formulated as follows.
\begin{prop}\label{prop1} (\cite{Lamperti}) If  $\{\bY(t), \ t\in \bbr \}$ is a strictly stationary stochastically continuous process and if for some $H>0$
$$
\bbx(t)=t^H \bY(\ln t), \quad t>0, \quad \bbx(0)=0,
$$
then $\{\bbx(t)\}$ is self-similar of order $H$. Conversely, every nontrivial self-similar process with $\bbx(0)=0$ is obtained in this way from a stationary process $$\bY(t)=e^{-H t}\bbx(e^t).$$
\end{prop}

In what follows we shall use the standard notation which is used, for example, in \cite{Samorod}, and we shall abbreviate the word self-similar by letters ss, so $H-ss$ process will stand for self-similar process with self-similarity (Hurst) index $H$. During last five decades the class of
$H-ss$ processes and some important subclasses such as $H-sssi$ processes (self-similar and with stationary increments) were studied in detail, list of reference would be long, we shall mention only two monographs, \cite {Samorod} and \cite{Embrechts1}, where one can find more references.

Self-similarity of processes is defined by transformations in time and space, and in the case of real valued processes of one-dimensional parameter $t$ these transformations are  simply scalings. Passing from one-dimensional parameter  to  multi-dimensional parameter $t$  one can see that there are several possibilities to generalize the notion of self-similarity. It seems that the most general definition of self-similarity is given in \cite{Kolodynski} for processes defined on a set$T$ and with values in $\bbr^m$. Let $G$ be a group of transformations of a set $T$ and let $C$ be a function defined on $G\times T$ in such a way that, for each $(g, t)\in G\times T$,  \ $C(g,t): \bbr^m \to \bbr^m$ is a bijection  satisfying the following relation
\begin{equation}\label{cocycle}
C(g_1g_2,t)=C(g_1,g_2(t))\circ C(g_2,t),
\end{equation}
for every $g_1, g_2 \in G, \  t\in T,$ and $C(e,t)=I_m$. Here $e$ stands for the unit element of $G$ and $I_m$ is the identity transformation of $\bbr^m$. Such function $C$ satisfying (\ref{cocycle}) is called a cocycle for the group action of $G$ on $T$.
\begin{definition}\label{groupss} (\cite{Kolodynski}) A  process $\bbx(t), \ t\in T,$ with values in $\bbr^m$ is called $G$-self-similar  with cocycle $C$ ($(G, C)$-ss), if, for each $g \in G$,
\begin{equation}\label{defgroup1}
\left \{ \bbx(g(t)), t\in T \right \}\stackrel{\rm d}{=}\left \{C(g,t) \bbx(t), t\in T \right \}.
\end{equation}
\end{definition}
In \cite{Kolodynski} there are several simple examples of $(G, C)$-ss processes with parameter set $T$ being $\bbr^d,$ $\bbr_+$, or $\bbn$, but the main aim of this paper was to characterize    $(G, C)$-ss strictly stable processes with values in $\bbr^m$ under general conditions on $T$ and $G.$
We shall consider mainly the case of $\bbr^m$-valued random fields defined on   $T=\bbr^d$ or $T=\bbr^d_+$, and cocycle independent of $t$: $C(g, \bt)=C(g).$  Also we shall consider  linear transformations in $\bbr^d$, taking $G\subseteq {\rm GL}(\bbr^d)$, then it is natural in the state space also to take linear bijections, that is, to require $C(g)\in {\rm GL}(\bbr^m)$. Here, as usual, ${\rm GL}(\bbr^d)$ stands for the linear group of invertible $d\times d$ matrices.
 Now we  reformulate Definition \ref{groupss}, taking cocycle  independent of $t$ and denoting it by letter $f$ instead of $C$.

\begin{definition}\label{rdgroupss} Let $G$ be some group of linear transformations of $\bbr^d$ and $f: G\to {\rm GL}(\bbr^m)$ be a function satisfying the relation $f(g_1g_2)=f(g_1)f(g_2)$, for all $g_1, g_2 \in G$. A $\bbr^m$--valued  process $\{\bbx(t), \ t\in \bbr^d,\}$ is called $(G, f)$-ss, if, for each $g \in G$,
\begin{equation}\label{defgroup2}
\left \{ \bbx(g(\bt)), \bt\in \bbr^d \right \}\stackrel{\rm d}{=}\left \{f(g) \bbx(\bt), \bt\in \bbr^d \right \}.
\end{equation}
\end{definition}
Let us denote by $\C$ the group of bijections in $\bbr^m$, generated by a cocycle $f$ and a group  $G$, i.e., $\C=\{f(g), g\in G\}\subseteq {\rm GL}(\bbr^m)$. It is natural to require that groups $G$ and $\C$ would be compatible in some sense. For example,
taking  $G=\{g={\rm {\rm diag}}(a_1, \dots, a_d), a_i>0, \ 1\le i \le d\}$ we shall require that $f(g)$ would be also diagonal, $f(g)={\rm diag}(f_1(\ba), \dots , f_m(\ba))$ (see Definition \ref{multiselfsim} below) and in this way in Proposition \ref{prop3}
 we obtain multi-self-similar random fields, which, in the case $m=1$, were introduced in  \cite{Genton}.
Taking $G=\{g=r^D, r>0\}$ with some fixed matrix $D$ on $\bbr^d$ with positive real parts of the eigenvalues, we obtain operator-scaling random field, in the case $m=1$, introduced in  \cite{Bierme} and in the general case $m\ge 1$ considered in  \cite{Li}.

Let us note, that, like in Lamperti Definition \ref{deflamperti}, in relations of type  (\ref{defgroup1}) or (\ref{defgroup2}) it is possible to add shift, and this possibility will be demonstrated bellow.

The goal of the paper is   to generalize Lamperti Theorems \ref{thm1}, \ref{thm2}, and Proposition \ref{prop1}, formulated above, to the more general setting of $\bbr^m$-valued random fields. We present several new results, but also we provide some earlier obtained results, trying to reflect all generalizations of the above mentioned  Lamperti results.
The rest of the paper  is divided into several sections. In section 2 we provide in some sense auxiliary results concerning regularly and slowly varying multivariate functions.  In section 3, taking different subgroups and cocycle functions in Definition \ref{rdgroupss} we describe several classes of self-similar random fields, while section 4 is devoted to domains of attraction of these random fields. In section 5  we consider the relation between the so-called phenomenon of the scale transition for random fields in the case $d=2$,  recently introduced in the papers \cite{Puplinskaite1} and  \cite{Puplinskaite2} and the domain of attraction of self-similar random fields.  In section 6 we discuss the Lamperti type transformation in general setting, and it is worth to note that only in this section we have the same level of generality as in Definition \ref{rdgroupss}.

\section{Regularly and slowly varying functions}

Regularly and slowly varying functions play an important role in problems connected with self-similarity of processes, therefore, considering random fields, we need some information about these functions of several variables. The notion of regular variation of real-valued and positive functions, defined on half-line, was introduced by J. Karamata in 1930 (see \cite{Karamata}). Now this theory is well developed, there are several monographs, devoted for this theory (see, for example, \cite{BINGHAM}). One can note that the regular variation of functions in many variables is investigated less than of univariate functions, and this can be explained by the fact that passing to multivariate functions we face with many possibilities even to define the notion of regular variation. Quite complete historical information one can find in \cite{Jakymiv}, see  remarks on p. 103, here we shall mention only two papers \cite{Karamata1} and \cite{Ostrogorski}. In the first above mentioned paper J. Karamata with his student initiated study of functions  $F: G\to \bbr_+$, where $G$ is a topological group with a given $G$-invariant filter $\U$ of open convex subsets of $G$ with countable basis and $\bbr_+$ is the multiplicative group of positive real numbers.
\begin{definition}\label{rvfunct-groups} (\cite{Karamata1}) Let $F: G\to \bbr_+$ be a continuous function defined above.  We say that this function  is regularly varying with respect to $\U$ if the following limit
\begin{equation}\label{regvar-groups}
\lim_{g\to \infty} \frac{ F(gh)}{F(g)}=\Phi(h),
\end{equation}
exists for every $h\in G$. Here  $g\to \infty$ means the convergence with respect to the filter $\U$.
\end{definition}
 From this definition  it is possible to prove in two lines the so-called Representation theorem: if (\ref{regvar-groups}) holds, then $\Phi$ is homomorphism of $G$ into $\bbr_+$, that is, $\Phi(h_1h_2)=\Phi(h_2)\Phi(h_1)$.

But the main stream of investigation in this area is connected with real-valued positive functions defined on cones in $\bbr^d$. We recall that a Borel set $\G \subset \bbr^d$ is called a cone with the vertex at a point $\ba\in \bbr^d$, if, for all $\bx\in \G,  \ \l \in \bbr, \l>0$,
$$
\ba+\l (\bx-\ba)\in \G.
$$
Let $\G$ be a cone with the vertex at $\B0$ and let $S=\G \setminus \{\B0\}$. In \cite{Jakymiv} a measurable function $f: S\to \bbr_+$  is called regularly varying at infinity with respect to $\G$ if for some fixed $\be\in S$ and all $\bx\in S$ the following limit exists
\begin{equation}\label{regvar-Jak}
\lim_{t\to \infty} \frac{ f(t\bx)}{f(t\be)}=\phi(\bx)>0.
\end{equation}
It is proved that the function $\phi$ from (\ref{regvar-Jak}) is homogeneous with the index $\rho \in \bbr$, that is, for $t>0$, $\phi(t\bx)=t^\rho \phi(\bx)$.

For the so-called homogeneous cones (to be defined below) T. Ostrogorski in \cite{Ostrogorski} had developed theory of regularly varying multivariate functions, based on ideas from \cite{Karamata1}. Let $V$ be an open convex cone in $\bbr^d$ and let $G$ be a subgroup of the general linear group of invertible matrices ${\rm GL}(\bbr^d)$ that leaves the cone $V$ invariant. The group $G$ is said to be transitive on $V$ if for every two points $\bx, \by \in V$ there is a $g\in G$ such that $\bx=g\by$, or, what is equivalent, if we  fix some $\be\in V$, then for every $\bx\in V$ there is $g\in G$ such that $\bx=g\be$. In this case we have a map $\pi: G\to V$, and if this map is one-to-one, the group $G$ is said to be simply transitive. A cone $V$ is called homogeneous if its group of linear automorphisms is transitive on $V$. The most important property of homogeneous cones, proved in \cite{Vinberg}, is that they always have a simply transitive group of automorphisms. On a homogeneous cone it is possible to define a product of its elements using group operation in $G$, also with any function $F: V\to \bbr_+$, using the map $\pi$, it is possible to relate  a function ${\bar F}: G\to \bbr_+$ by putting ${\bar F}=F\circ \pi$. Then regular variation of $F$ is defined via regular variation of ${\bar F}$ using Definition \ref{rvfunct-groups} (with some filter $\U$). We shall not provide here  the strict formulations of results in this direction, since it would require many new notions from algebra (such as Lie algebra of a group, exponential mapping , rank of a cone, etc.), instead of that we can formulate important message from \cite{Ostrogorski}: taking different cones, different simple transitive groups of automorphisms and filters, we may get different classes of regularly varying multi-variate functions.

In the context of generalization of Lamperti theorems we are interes-ted in the particular cone $\bbr^d_+=(0, \infty)^d$, which is clearly homogeneous cone. As $G$ we can take diagonal matrices with positive entries, that is $G=\{{\rm diag} (t_1, \dots , t_d), \ t_i>0\}$, and it is easy to verify that $G$ is simply transitive on $\bbr^d_+$. Moreover, $G$ is homomorphic to $\bbr^d_+$: if $g={\rm diag} (t_1, \dots , t_d)\in G$, then $\pi (g)= (t_1, \dots , t_d)\in \bbr^d_+$, therefore in this cone product of two elements can be defined coordinate-wise. Thus we arrive at the following definition and proposition.
\begin{definition}\label{rvfunctions}   We say that  $f: \bbr_+^d\to \bbr_+$  is a coordinate-wise regularly varying  function (c.r.v.f.) if, for all $\bx\in \bbr_+^d$,
\begin{equation}\label{regvar}
\lim_{\bt\to \infty} \frac{ f(\bt\cdot \bx)}{f(\bt)}=l(\bx),
\end{equation}
with some continuous  function $l$. Here and in what follows $\bt\to \infty$ means $\min_{1\le i\le d} t_i\to \infty$.

 We say that  $L: \bbr_+^d\to \bbr$  is coordinate-wise slowly varying  function (c.s.v.f.) if, for all $\bx\in \bbr_+^d$,
\begin{equation}\label{slowvar}
\lim_{\bt\to \infty} \frac{ L(\bt\cdot \bx)}{L(\bt)}=1.
\end{equation}
\end{definition}

\begin{prop}\label{prop2} Let $f: \bbr_+^d\to \bbr_+$ be c.r.v.f. . Then there exist positive $H_i, 1\le i\le d,$ such that
\begin{equation}\label{regvar1}
f(\bt)=\prod_{i=1}^dt_i^{H_i}L(\bt),
\end{equation}
where $L$ is c.s.v.f. .
\end{prop}

{\it Proof. }
Although this proposition can be derived from general results from \cite{Ostrogorski}, we shall give very simple elementary proof without any use of sophisticated notions from algebra.
Taking $\ba, \bb\in \bbr_+^d$ and  passing to the limit as $\bt\to \infty$ in the identity
$$
 \frac{ f(\ba \cdot \bb \cdot \bt)}{f(\bt)}= \frac{ f(\ba \cdot \bb \cdot \bt)}{f( \bb \cdot \bt)} \frac{ f( \bb \cdot \bt)}{f(\bt)}
$$
we shall get that for all $\ba, \bb\in \bbr_+^d$
\begin{equation}\label{identity1}
l(\ba \cdot \bb )=l(\ba)l(\bb).
\end{equation}
Taking $a_i=b_i=1, i=2, \dots, d$  we get
$$
 l(a_1b_1, 1, \dots, 1)=l(a_1, 1, \dots, 1)l(b_1, 1, \dots, 1)
$$
From this functional equation, taking into account that function $l( \cdot, 1, \dots, 1)$ is continuous we conclude that $ l(a_1, 1, \dots, 1)= a_1^{H_1}$.
From (\ref{identity1}) we obtain
\begin{equation}\label{identity2}
l(a_1, b_2, \dots , b_d )=l(a_1, 1, \dots, 1)l(1,  b_2, \dots , b_d )=a_1^{H_1}l(1,  b_2, \dots , b_d ).
\end{equation}
Let us denote $l_1(a_2, \dots , a_d):=l(1,  a_2, \dots , a_d )$. Repeating the same procedure with function $l_1$ for the  variable $a_2$ we shall get
\begin{equation}\label{identity3}
l_1(a_2, a_3, \dots , a_d )==a_2^{H_2}l_1(1,  a_3, \dots , a_d ).
\end{equation}
From equalities (\ref{identity2}) and (\ref{identity3}) we get
$$
l(a_1, a_2, \dots , a_d )==a_1^{H_1}a_2^{H_2}l(1, 1,  a_3, \dots , a_d ),
$$
therefore, continuing in this way, we shall arrive at the equality
$$
l(a_1, \dots , a_d )==\prod_{i=1}^d a_i^{H_i}.
$$

Let $f$ and $g$ be two functions satisfying (\ref{regvar}) with the same function $l$. Denote
$$
L(\bt)=\frac{f(\bt)}{g(\bt)}.
$$
From (\ref{regvar}) it follows that
$$
\lim_{\bt\to \infty} \frac{ f(\ba \cdot \bt)g(\bt)}{g(\ba \cdot \bt)f(\bt)}=1,
$$
but this is the relation (\ref{slowvar}) for the function $L$, i.e., function $L$ is c.s.v.. Since  the function $g(\bt)=\prod_{i=1}^d t_i^{H_i}$ clearly satisfies (\ref{regvar}), whence it follows that any function, satisfying (\ref{regvar}) is of the form (\ref{regvar1}).
\halmos

If we would like to apply the same approach from \cite{Ostrogorski} to the cone $\bbr_+^d$ and taking as the definition of regularly varying function  the relation  (\ref{regvar-Jak}) (as it is done in \cite{Jakymiv}) instead of (\ref{regvar}), we shall fail, since in this case the group of transformations $G$ would be diagonal matrices with equal positive entries, and it is easy to see that such group is not transitive on the cone  $\bbr_+^d$. Clearly, relation (\ref{regvar-Jak}) reflects the radial behavior of the function $f$,  since we compare the growth of the function on a ray $\{ t\bx, \ t>0\}$, for a fixed $\bx$ with the growth on the other fixed ray
$\{ t\be, \ t>0\}$.   Therefore we shall get a class of functions, different from functions from (\ref{regvar1}). It will be convenient to fix $\be =(d^{-1/2}, \dots , d^{-1/2})$. For $\bx \in \bbr^d$, let us denote $r_\bx=||\bx||, \ a_\bx=\bx/||\bx||$.
\begin{definition}\label{radRVF}   We say that  $f: \bbr_+^d\to \bbr_+$  is a radially regularly varying  function (r.r.v.f.) if, for all $\bx\in \bbr_+^d$, the relation (\ref{regvar-Jak}) holds.

We say that $L: \bbr_+^d\to \bbr$  is radially slowly varying  function (r.s.v.f.) if, for all $\bx\in \bbr_+^d$,
\begin{equation}\label{radslowvar}
\lim_{t\to \infty} \frac{ L(t \bx)}{L(t\be)}=\l(a_\bx)
\end{equation}
with some non-negative function $\l$, defined on the part of the unit sphere $\{\bx \in \bbr^d_+: ||\bx||=1\}$ and satisfying $\l(a_{\be})=1$.
\end{definition}

  In \cite{Jakymiv} more narrow class of slowly varying functions on  cones is defined, since it is required the relation  (\ref{radslowvar}) with $\l(a_\bx)\equiv 1$.  Also it is possible to think about more general,  comparing with Definition \ref{radRVF}, r.r.v.f., when the exponent of regular variation also depend on the direction, that is, it is possible to consider functions of the form
$$
f(\bx)=r_\bx^{\rho(a_\bx)}L(\bx),
$$
where $L$ is a r.s.v.f.. Using the same idea (that on each ray the behavior can be different) we can consider more general r.s.v.f.  of the form $L(\bx)=l_{a_\bx}(r_\bx),$ where for each $a_\bx$  on the unit sphere $l_{a_\bx}$ is a univariate s.v.f. . But at present we do not know in which problems such general functions can appear.

\begin{prop}\label{prop10} Let $f: \bbr_+^d\to \bbr_+$ be r.r.v.f. . Then there exist  $\rho \in \bbr,$ such that
\begin{equation}\label{radregvar1}
f(\bx)=r_\bx^{\rho}L(\bx),
\end{equation}
where $L$ is a r.s.v.f. . $L$ is a radially slowly varying function, if and only if
 \begin{equation}\label{radSVF}
L(\bx)=K(\bx)l(r_\bx)\l(a_\bx),
\end{equation}
with some univariate s.v.f. $l$,  some non-negative function $\l$, defined on the part of unit sphere $\{\bx \in \bbr^d_+: ||\bx||=1\}$, and some positive function $K$  having the expression (\ref{slowvarform-Jak}).
\end{prop}

{\it Proof. } Let us prove at first (\ref{radSVF}). Let $L$ be a function satisfying (\ref{radslowvar}) and let us take the function
$$
L_1(\bx)= l(r_\bx)\l (a_\bx),
$$
with some univariate s.v.f. $l$ and the same as in (\ref{radslowvar}) function $\l$. It is easy to verify that $L_1$ satisfy (\ref{radslowvar}). Let us denote $K(\bx)=L(\bx)/L_1(\bx)$. Taking into account that $r_\be =1$ and $\l (a_\be)=1$ it is not difficult to verify that, for all $\bx \in \bbr_+^d$,
$$
\lim_{t\to \infty} \frac{ K(t\bx)}{K(t\be)}=1.
$$
Thus, we get that any function satisfying (\ref{radslowvar}) is of the form (\ref{radSVF})
with a function $K$, satisfying (\ref{radslowvar}) with $\l(a_\bx)\equiv 1$. As it was mentioned after formulation of Definition \ref{radRVF}, this is exactly the definition of the s.v.f. given in \cite{Jakymiv}, where it is given the general form of such functions (see Theorem1.1.4 in \cite{Jakymiv}):
\begin{equation}\label{slowvarform-Jak}
K(\bx)=\exp \left (\eta (\bx)+\int_a^{r_{\bx}} \frac{\e (u)}{u} du\right ),
\end{equation}
with some $a>0, $ and some functions $\eta (\bx)\to C, \ C\in \bbr,$ and $\e(t)$, for any $k\ge 0$, satisfying $\e (t)=o(t^{-k})$. Thus, we proved  (\ref{radSVF}).

Now it is easy to prove (\ref{radregvar1}).
 In \cite{Jakymiv} it is proved that if a function $f$ satisfies (\ref{regvar-Jak}), then  $\phi $ is homogeneous, this means that $\phi (\bx)=r_\bx^\rho \phi (a_\bx)$ with some $\rho \in \bbr$.
Let us take the function
$$
f_1(\bx)=r_\bx^\rho l(r_\bx)\phi (a_\bx),
$$
with some univariate s.v.f. $l$ and with the same as in (\ref{regvar-Jak}) function $\phi$, and denote $L(\bx)=f(\bx)/f_1(\bx)$. As in the proof above, taking into account that $r_\be =1$ and $\phi (a_\be)=1$, one can get  that, for all $\bx \in \bbr_+^d$,
$$
\lim_{t\to \infty} \frac{ L(t\bx)}{L(t\be)}=1.
$$
Therefore, we get that any function satisfying (\ref{regvar-Jak}) is of the form
\begin{equation}\label{regvarform}
f(\bx)=r_\bx^\rho l(r_\bx)\phi (a_\bx)L(\bx),
\end{equation}
with a function $L$, satisfying (\ref{radslowvar}) with $\l(a_\bx)\equiv 1$. But we know that this function is of the form given in (\ref{slowvarform-Jak}).
 Combining this representation (\ref{slowvarform-Jak}), (\ref{regvarform}), and (\ref{radSVF}) we obtain (\ref{radregvar1}).
\halmos

In Definitions \ref{rvfunctions} and \ref{radRVF} we had introduced two classes of  regularly varying functions and two classes of slowly varying functions. We denote functions, satisfying (\ref{slowvar})  by $\cl_1$ and by $\cl_2$ we denote functions, satisfying (\ref{radslowvar}). From definitions it follows that $\cl_1 \subset \cl_2$.  Regularly varying functions from both introduced classes have the same structure
$$
f(\bx)=g(\bx)L(\bx),
$$
where $g(\bx)=\prod_{i=1}^dx_i^{H_i}$ in the case of c.r.v.f and $g(\bx)=||\bx||^\rho$ in the case of r.r.v.f., while $L$ is slowly varying function from the corresponding class. Thus, for the complete description of both classes of regularly varying functions, it remains to have such description for the class $\cl_1$ , but it seems that at present it is not available.
We can provide some examples of c.s.v.f. and r.s.v.f., constructed from univariate slowly varying functions. Let $l_i(x), 1\le i\le d$ be arbitrary univariate s.v.f. and let $\phi$ be some function defined on the unit sphere of $\bbr^d$. Denote
$$
L_1(\bx)=\prod_{i=1}^d l_i(x_i), \quad L_2(\bx)=\sum_{i=1}^d l_i(x_i), \quad L_3(\bx)=l_1(r_\bx)\phi (a_\bx),
$$
where we recall $r_\bx= ||\bx||, a_\bx=\bx/||\bx||$. An easy exercise in calculus shows that $L_i \in \cl_1$, for $i=1,2$, and $L_3 \in \cl_2$.

Different approach to define multivariate  regularly varying functions is used in papers \cite{Meerschaert} and \cite{Meerschaert1}. At first  the notion of regular variation of a function $F:\bbr_+ \to {\rm GL}(\bbr^d)$ is introduced: we say that $F$ is regularly varying at infinity with index $D$ if, for all $\l>0$,
\begin{equation}\label{opregvar}
\lim_{r\to \infty} F(\l r)F(r)^{-1}=\l^D,
\end{equation}
where $D$ is some (possible singular) linear operator on $\bbr^d$ and convergence of matrices is in the usual norm in  ${\rm GL}(\bbr^d)$. Here $\l^D$ stands for the operator $\exp (\ln \l \cdot D)$ and $\exp (A)=\sum_{n=0}^\infty A^n/n!$. As an example of regularly varying function one can take a function $F(r)=r^E$ with some  matrix $E$. Function $L:\bbr_+ \to {\rm GL}(\bbr^d)$  is slowly varying at infinity if, for all $\l>0$,
\begin{equation}\label{opslowvar}
\lim_{r\to \infty} L(\l r)L(r)^{-1}=I,
\end{equation}
 where $I$ is the identity matrix in ${\rm GL}(\bbr^d)$. It seems that it is not known the form of general slowly varying function with values in ${\rm GL}(\bbr^d)$, we can only provide the following simple example of such function. Let $l_i(r), 1\le i\le d$, be arbitrary univariate s.v.f.'s, then it is easy to see that the function $L(r)= {\rm diag} (l_1(r), \dots , l_d(r))$ is slowly varying function in the sense of (\ref{opslowvar}). But there is one essential difference between regularly varying functions with values in ${\rm GL}(\bbr^d)$ and functions defined in Definitions \ref{rvfunctions} and \ref{radRVF}.
 Due to the non-commutativity of  ${\rm GL}(\bbr^d)$ it is not possible to get the following representation theorem: any function, regularly varying at infinity with index $D$, is of the form $F(r)=r^D L(r)$ where $L:\bbr_+ \to {\rm GL}(\bbr^d)$ is some slowly varying function. Also for the same reason it is no longer true that product of two  regularly varying functions is again regularly varying function, while for regularly varying functions, defined in Definitions \ref{rvfunctions} and \ref{radRVF}, such property holds.

  Having defined regularly varying functions with values in ${\rm GL}(\bbr^d)$, the regularly varying function $f:\bbr^d\setminus \{0\}\to \bbr_+$ is defined in \cite{Meerschaert} using the regularly varying function $F$  with index $D$ and univariate regularly varying function $l$ with index of regularity $\b$ in the following way. It is said that
$f:\bbr^d\setminus \{0\}\to \bbr_+$ is regularly varying if there exist regularly varying functions $F:R_+ \to {\rm GL}(\bbr^d)$ and $l: \bbr_+ \to \bbr_+$ such that
\begin{equation}\label{opregvarrd}
\lim_{r\to \infty} \frac{f(F(r)^{-1}\bx_r)}{l(r)}=\phi(\bx)>0,
\end{equation}
for any $\bx_r\to \bx$ in $\bbr^d\setminus \{0\}$. If all eigenvalues of $D$ have positive real parts, then $|| F(r)^{-1}\bx_r||\to 0$ if $r\to \infty$ and   $\bx_r\to \bx$, and relation (\ref{opregvarrd}) defines regular variation of $f$ at $0$, while, in the case of all negative real parts of eigenvalues of $D$, this relation defines the behavior of $f$ at infinity.


\section{Self-similar random fields}

As it was mentioned above, taking particular groups $G \subset {\rm GL}(\bbr^d)$ and cocycle functions $f$ in Definition \ref{rdgroupss}, we get different classes of self-similar random fields. We start with simple case of diagonal matrices for transformations both in time and space, and since the group of diagonal matrices with positive entries  is homomorphic to $\bbr^d_+$, we shall use the same notation as in the previous section: if   $g={\rm diag} (a_1, \dots , a_d)$, with $\ba:=(a_1, \dots , a_d)>\B0$,   then instead of $g(\bt)$ we shall write $\ba\cdot \bt$. In the following definition and theorem we can consider not only $\bbr^d$, but also $\bbr^d_+$, since for the diagonal transformation $g$ with positive entries we have $g(\bt)\in \bbr^d_+$ if $\bt \in \bbr^d_+$. Therefore, in these two statements $T$ stands for $\bbr^d$ or $\bbr^d_+$.
\begin{definition} \label{multiselfsim} A  random field $\{\bbx(\bt)=(X_1(\bt), \dots , X_m(\bt)), \ \bt \in T\}$ with values in $\bbr^m$, $m\ge 1$, is called multi-self-similar (m.s.s.), if it is stochastically continuous and for any $\ba>\B0$ there exist  function $\bbf(\ba)= {\rm diag} (f_1(\ba), \dots , f_m(\ba))$  such that
 \begin{equation}\label{defmss}
 \{\bbx(\ba\cdot \bt),  \ \bt \in T\}\stackrel{\rm d}{=} \{\bbf(\ba)\cdot\bbx(\bt),  \ \bt \in T\}.
  \end{equation}
A random field $\bbx$  is called wide-sense multi-self-similar (w.m.s.s.), if it is stochastically continuous and there exist functions $\bbf(\ba)$ and $\bg(\ba)=(g_1(\ba), \dots , g_m(\ba))$  such that
 \begin{equation}\label{widesmss}
 \{\bbx(\ba\cdot \bt),  \ \bt \in T\}\stackrel{\rm d}{=} \{\bbf(\ba)\cdot\bbx(\bt)+\bg(\ba),  \ \bt \in T\}.
 \end{equation}
\end{definition}
\begin{teo}\label{prop3} If a random field $\{\bbx(\bt),  \ \bt \in T\}$ is  non-degenerate and m.s.s., then there exist  vectors  $\bbh^{(j)}=(H_1^{(j)}, \dots, H_d^{(j)}), \ j=1, \dots , m,$  with non-negative coordinates such that
 \begin{equation}\label{fexpres}
f_j(\ba)=\prod_{i=1}^d a_i^{H_i^{(j)}}, \ j=1, \dots , m,
\end{equation}
and $X(\B0)=0$. If a random field is w.m.s.s., then $\bbf$ has the same expression (\ref{fexpres}) and
\begin{equation}\label{gexpres}
g_j(\ba)=D_j\left (1-f_j(\ba) \right ),
\end{equation}
where $D_j=X_j(\B0)$.
If $H_i^{(j)}=0 $ for all $i$, then $g_j(\ba)\equiv 0$ and $X(\bt)\equiv X(\B0)$ a.s..

If $H_i^{(j)}=0$ for some collection of indices $A_j\subset \{1, \dots , d\}$, then denoting $Y_j(t_i,i\in A)= X_j(\bt)$ we have  $Y_j(t_i,i\in A_j)\equiv Y_j(0,\dots ,0)$, i.e. the process $X_j$ is constant with respect to those variables $t_i$ for which  $H_i^{(j)}=0$.
\end{teo}

If we denote by letter ${\bar \bbh}$ the $m\times d$ matrix formed by vectors $\bbh^{(j)}, j=1, \dots , m,$ then  the notations ${\bar \bbh}$-m.s.s. and ${\bar \bbh}$-w.m.s.s.  will stand  for m.s.s. and w.m.s.s.  random fields, respectively.
As it was mentioned above, m.s.s. random fields in the case $m=1$ and for $T=\bbr^d$ were introduced in \cite{Genton}, but even in the case  $m=1$ there is a difference between our Definition \ref{multiselfsim} and the Definition 1 in \cite{Genton}: in the latter definition particular form of the function $f_1$  is required, while we, like in the original Lamperti paper, prove that  the function $f_1$ is of that particular form, also we introduce more general definition allowing shift.

{\it Proof of Theorem \ref{prop3}}.  As in the proof of Lamperti theorem \ref{thm1}, at first, using the stochastic continuity of the process $\bbx$, we can show that functions $\bbf$ and $\bg$ are continuous.

 Assuming that the process $\bbx$ is w.m.s.s., from definition (\ref{widesmss}) we can write the following two equalities
 $$
 \PP \{\bbx(\ba \cdot \bb \cdot \Bl1)\le \by\}=\PP \left \{\bbx(\Bl1)\le (\bbf(\ba \cdot \bb))^{-1} (\by-\bg(\ba \cdot \bb))  \right \}
 $$
and
$$
\PP \{\bbx(\ba \cdot \bb)\le \by\}=\PP \left \{\bbx(\bb)\le (\bbf(\ba))^{-1} \left (\by-\bg(\ba)\right )  \right \}
$$
$$
=\PP \left \{\bbx(\Bl1)\le (\bbf(\bb))^{-1}\left ((\bbf(\ba))^{-1} \left (\by-\bg(\ba)\right )-\bg(\bb)\right ) \right \}.
$$
Here $(\bbf(\ba))^{-1}$ is the inverse matrix of the diagonal matrix, therefore it is itself diagonal, $(\bbf(\ba))^{-1}= {\rm diag} (f_1^{-1}(\ba), \dots , f_m^{-1}(\ba))$.
 From these two relations, equating vectors $(\bbf(\ba \cdot \bb))^{-1} (\by-\bg(\ba \cdot \bb))$ and $(\bbf(\bb))^{-1}\left ((\bbf(\ba))^{-1} \left (\by-\bg(\ba)\right )-\bg(\bb)\right )$, we conclude
\begin{equation}\label{widesmss1}
 f_j(\ba \cdot \bb)=f_j(\ba)f_j(\bb), \ j=1, \dots , m,
 \end{equation}
 \begin{equation}\label{widesmss2}
g_j(\ba \cdot \bb)=g_j(\ba) +g_j(\bb) f_j(\ba), \ j=1, \dots , m.
\end{equation}
We see that all equations (\ref{widesmss1}) and (\ref{widesmss2}) with respect to $j$ have the same form, so it is sufficient to solve only one for some fixed  $1\le j\le m$.
Clearly, if the process $\bbx$ is m.s.s., then we have only  the equations (\ref{widesmss1}). These equations (for any $j$) are exactly of the form (\ref{identity1}), therefore, $f_j(\ba)=\prod_{i=1}^d a_i^{H_i^{(j)}},$ and we have proved (\ref{fexpres}). Now we must solve the equation
 \begin{equation}\label{widesmss3}
g_j(\ba \cdot \bb)=g_j(\ba) +g_j(\bb)\prod_{i=1}^d a_i^{H_i^{(j)}}.
\end{equation}
It is not difficult  to verify  that the function $g_j(\ba)= D_j(1-\prod_{i=1}^d a_i^{H_i^{(j)}}) $ with some constant $D_j$ satisfies the equation (\ref{widesmss3}). From relation (\ref{widesmss}), taking $\bt=\B0$ we get
$$
  X_j(\B0)\stackrel{\rm d}{=}\prod_{i=1}^d a_i^{H_i^{(j)}} X_j(\B0)+ g_j(\ba),
$$
whence we obtain that $D_j=X_j(\B0)$. The formal proof that only this function $g_j$ is the solution of (\ref{widesmss3}) can be carried along the same lines as  the proof in \cite{Lamperti} in the case $d=1$. Namely, making the change of variables $b_i=e^{u_i}, a_i=e^{v_i}, i=1, \dots , d, $ and denoting $d_j(\bu)=g_j(e^{\bu})$, where $e^{\bu}=\left (e^{u_1}, \dots, e^{u_d} \right )$, we must solve the equation
\begin{equation}\label{widesmss4}
d_j(\bu+\bv)=d_j(\bv) +d_j(\bu)\exp \{\sum_{i=1}^d v_i H_i^{(j)} \}.
\end{equation}
We must prove that the function
\begin{equation}\label{solution}
d_j(\bu)= D_j\left (1-\exp \left \{\sum_{i=1}^d u_i H_i^{(j)} \right \}\right ), \quad \bu\in \bbr^d,
\end{equation}
with some $D_j$ is the solution of (\ref{widesmss4}). In the proof of this fact we use the same method which was used in the original Lamperti proof in \cite{Lamperti}. Namely, we take a vector $\bv$ with coordinates  from the set $\{0, 1\}$ and we write a collection of equalities obtained from (\ref{widesmss4}). Then, using induction, we prove that (\ref{solution}) is a solution for $\bu$ with integer coordinates, then in a similar way we prove for rational coordinates, and it remains to use the continuity of $d_1$. We have proved  (\ref{gexpres}).

It remains to consider the cases where some of coordinates of vectors $\bbh^{(j)}, j=1, \dots , m,$ are equal to $0$, and, again, it is sufficient to consider only one vector $\bbh^{(j)}.$

 The case  $H_i^{(j)}=0 $ for all $i$ is easy, since in this case the relation (\ref{widesmss}) for the $j$th coordinate becomes
 $$
 \{X_j(\ba \cdot \bt)\}\stackrel{\rm d}{=} \{X_j(\bt)+g_j(\ba)\}
 $$
 and taking $\bt=\B0$ we get $g_j(\ba)\equiv \B0$. Then taking $\ba=\B0$ we get $X_j(\bt)\equiv X_j(\B0)$ a.s.

For simplicity of writing let us consider the case where $A_j=\{1\}$, i.e.
we consider the case $H_1^{(j)}=0, \ H_i^{(j)}>0, \ i\ge 2$, then we have
$$
g_j(\ba)=D_j\left (1-\prod_{i=2}^d a_i^{H_i^{(j)}}\right )
$$
and taking $\ba=(a_1, 1, \dots , 1)$ from relation (\ref{widesmss}) we get
$$
\{X_j(a_1t_1, t_2, \dots , t_d)\}\stackrel{\rm d}{=} \{X_j(t_1, t_2, \dots , t_d)\},
$$
since $g_j(\ba)=0$ and $f_j(\ba)=1$ for this $\ba$. In the last relation taking $a_1=0$  we get $X_1(t_1, t_2, \dots , t_d)=X(0, t_2, \dots , t_d)$ for all $\bt$.

\halmos

Taking another group of transformations of "time" parameter $\bt$  we get   multivariate operator self-similar random fields, introduced in \cite{Bierme}  and \cite{Li}. Instead of the group of diagonal matrices with positive entries, which is homomorphic to $\bbr^d_+$, one-parameter group $G=\{r^E, \ r>0,\}$ with some fixed $d\times d$ matrix $E$ is used in definition of this new class of self-similar fields. Now we can consider only the case of $T=\bbr^d$, since the requirement $r^E\bt\in \bbr^d_+$, if $\bt \in \bbr^d_+$, puts too strong requirement for the matrix $E$ (only diagonal matrices with positive entries on the diagonal). Let us denote by $Q(\bbr^d)$ the set of invertible $d\times d$ matrices whose eigenvalues have positive real parts and by $M(\bbr^d)$ the set of $d\times d$ matrices $A$ such that all eigenvalues of $A$ have nonnegative real parts and every eigenvalue of $A$ with real part equal to zero (if it exists) is a simple root of the minimal polynomial of $A$.  The following definition is given in \cite{Li}.
 \begin{definition} \label{selfsimop} Let $E\in Q(\bbr^d)$. A $\bbr^m$-valued random field $\{\bbx(\bt),  \ \bt\in \bbr^d\}$ is called wide-sense operator self-similar (w.o.s.s. ) with time-variable scaling exponents $E$, if for any $r>0$ there exist an $m\times m$ matrix $B(r)$ (which is called a state space scaling operator)  and a function $\ba_r(\cdot):\bbr^d \to \bbr^m$ such that
\begin{equation}\label{defgroup4}
\left \{ \bbx(r^E\bt)), \bt\in \bbr^d \right \}\stackrel{\rm d}{=}\left \{B(r) \bbx(\bt)+\ba_r(\bt), \bt\in \bbr^d \right \}.
\end{equation}
If, in addition, $\ba_r(\bt)\equiv 0$, then $\bbx$ is called  operator self-similar (o.s.s. ) with time-variable scaling exponents $E$. If the function $\ba_r$ in (\ref{defgroup4}) does not depend on $\bt$, i.e., a process $\bbx$ satisfies
$$
\left \{ \bbx(r^E\bt)), \bt\in \bbr^d \right \}\stackrel{\rm d}{=}\left \{B(r) \bbx(\bt)+\ba(r), \bt\in \bbr^d \right \},
$$
the process $\bbx$ is called  operator self-similar in the sense of Hudson and Mason (HMo.s.s. ).
\end{definition}
In \cite{Li}  descriptions of functions $B(r), a_r(\bt), a(r)$, appearing in Definition \ref{selfsimop},  are given. We formulate one such result.

\begin{teo}\label{thm3} (\cite{Li}) Let $\{\bbx(\bt), \bt \in \bbr^d\}$ be a stochastically continuous and proper $\bbr^m$-valued w.o.s.s.  random field with time-variable scaling exponent $E\in Q(\bbr^d)$. There exist a matrix $D\in M(\bbr^m)$ and a function $b_r(\bt): (0, \infty)\times\bbr^d \to \bbr^m$ which is continuous at every $(r, \bt)\in (0, \infty)\times\bbr^d $ such that for all constants $r>0$
\begin{equation}\label{defgroup6}
\left \{ \bbx(r^E\bt)), \bt\in \bbr^d \right \}\stackrel{\rm d}{=}\left \{r^D \bbx(\bt)+b_r(\bt), \bt\in \bbr^d \right \}.
\end{equation}
Furthermore, $\bbx(\B0)=\ba$ a.s. for some constant vector $\ba \in \bbr^m$ if and only if $D\in Q(\bbr^m)$. In this latter case, we define $b_0(\bt)\equiv \ba$ for all  $\bt \in \bbr^d$, then the function $(r, \bt) \to b_r(\bt)$ is continuous on $[0, \infty)\times\bbr^d. $
\end{teo}

The matrix $D$ in the representation (\ref{defgroup6}) is called space-scaling exponent. Similar results are formulated for o.s.s.  and HMo.s.s.  random fields.  For all these random fields the main characteristics are these two exponents $E$ and $D$ and shortly we shall denote  an o.s.s.  random field with exponents $E$ and $D$ by $(E, D)$-o.s.s.  random field. These three classes of  operator self-similar random fields are more complicated comparing with multi-self-similar random fields. One complication steams from the fact that for a given exponent $E$ exponent $D$ may be not unique, and in a recent paper \cite{Didier} this problem is considered.

 Let us  compare  Definition \ref{selfsimop} of an o.s.s.  random field with Definition \ref{multiselfsim} (with $T= \bbr^d$ ) of ${\bar \bbh}$-m.s.s.   random field, considered in Theorem \ref{prop3}. 
 Suppose that $\bbx$ is  ${\bar \bbh}$-m.s.s.   random field, thus, taking into account Theorem \ref{prop3}, it satisfies (\ref{defmss}) with $\bbf(\ba)$ from (\ref{fexpres}). For any $r>0$ and arbitrary fixed $E={\rm diag} (e_1, \dots , e_d),$ let us take $\ba =r^E$. From (\ref{defmss}) and (\ref{fexpres})    we see that ${\bar \bbh}$-m.s.s. random field  will be $(E, D)$-o.s.s.   random field with a given diagonal matrix $E$ and $D$ being diagonal with elements $b_j=\sum_{i=1}^d e_iH_i^{j}, j=1, \dots, m$, on diagonal.

\smallskip

  But even  if we  restrict class of matrices in definition of operator-scaling random fields (\ref{defgroup4}) by diagonal matrices $E={\rm diag} (e_1, \dots , e_d)$ with positive entries on diagonal, still we have difference between these two definitions. The difference is between the convergence $ \ba \to \infty$ (which is $\min (a_i) \to \infty$), while $r^D \to \infty$ means $r\to \infty$, thus the paths $\ba \cdot \bt$ and $r^D\bt$ of these two kinds of convergence  in $\bbr^d$ are different.

\section{Domains of attraction of self-similar random fields}

The problem of generalizing Lamperti Theorem \ref{thm2} for random fields is quite difficult. Having general Definitions \ref{groupss} and \ref{rdgroupss} of self-similar random field, defined on an abstract set $T$ by means of some group of transformations of $T$, one would like to define domain of attraction, having the same generality. One can try to define domain of attraction of a process $\{\bbx(t), t\in T\}$ with values in $\bbr^m$ in the following way.
\begin{definition} \label{gendomattr}
Let $G$ be some topological group of transformations of $T$ with some filter $\U$. Let us say that a $\bbr^m$-valued process $\{\bY(t), t\in T\}$ belongs to the domain of attraction, defined by  $G$, of the process $\{\bbx(t), t\in T\}$, if there exist a functions $F: G\to {\rm GL}(\bbr^m)$ and $\ba: G\to \bbr^m$ (with some properties) such that
\begin{equation}\label{lim2}
\lim_{g \to \infty}\left \{F(g)\left (\bY(g(t))+\ba(g)\right ), t\in T\right \} \fdd \{ \bbx(t), t\in T\},
\end{equation}
where  $g\to \infty$ means the convergence with respect to the filter $\U$.
\end{definition}
 Then, from the relation (\ref{lim2}), one can try to prove that the process $\bbx$ is $(G, f)$-ss, with some cocycle $f$ in the sense of Definition \ref{groupss}, and to establish the relation between $F$ and $f$. Let us note that in this definition centering (the term $\ba(g)$), like in the original definition of Lamperti (see (\ref{lim1})) or later result in Theorem 5 from \cite{Hudson},  is independent of $t \in T$.
  But  simple examples shows that such assumption is unnatural and the centering depends on $t$. Namely, let us consider a sequence of i.i.d. random variables $X_1, \dots , X_n, \dots$ belonging to the domain of attraction of a stable distribution with index $1<\a \le 2$ and with $EX_1=a\ne 0$. Then it is clear that the centering for the process $S_n(t)=\sum_{i=1}^{[nt]}X_i$ is dependent on $t$ and is $[nt]a$ and the limit process is $H$-ss with $H=1/\a$. The same example shows that there is a misprint in the expression of the function $\bk$ from  (\ref{lim1}): instead of $\bk(a)=\omega(a)a ^H L(a)$ there should be $\bk(a)=\omega(a)a  L(a)$. That this is misprint confirms the following statement in the proof (see p. 73 of \cite{Lamperti}): ``...and $\bk(a)^H$ to be $\omega a^H$".
  More generally, if we have a process $Y(t)$  with finite mean, then it is clear that as the centering must be taken $EY(t)$ and we can use the centered process ${\bar Y}(t)=Y(t)-EY(t)$. Therefore, we have two possibilities to modify (\ref{lim2}): to use $\ba (g, t)$ instead of $\ba (g)$ or to set $\ba (g)= \B0$. In the paper we choose the second possibility and define domains of attraction without centering.

   Another remark, concerning Definition \ref{gendomattr}, is that such formulation of the problem is too general and some additional assumptions on a group of transformations $G$ must be adopted in order to avoid some trivial situations, like the following one. If we take $T=\bbr$ and the group of translations as $G$ and the cocycle function equal to 1, then self-similarity with respect to this group is nothing else as stationarity. It is evident that the domain of attraction of a stationary process $X(t), t\in \bbr,$ defined as above,  is trivial and
consists of the process $X$ itself. Therefore, we restrict our investigation by taking $T=\bbr^d$ or $T=\bbr^d_+$ and two particular   groups $G$, which were considered in Section 3. In each case  we give a separate definition of the domain of attraction. In  Definition  \ref{gendomattr} the domain of attraction depends on $G$, but there is no requirement for the function $F$, now we put some requirements on this function, and these requirements, naturally,  are different for different groups of transformations $G$. Such approach can be justified by  the following considerations. It is well-known that
domains of attraction are closely related with summation theory.  Suppose that $\xi = \{\xi_{\bi}, \bi\in \bz^d\}$ is a $\bbr^m$-valued stationary random field and
\begin{equation}\label{sumfields}
\bbs_{[\bn \cdot \bt]}(\xi)=\sum_{\bi=\Bl1}^{[\bn \cdot \bt]}\xi_{\bi}, \quad   [\bn \cdot \bt]=([n_1 t_1], \dots , [n_dt_d]), \   \sum_{\bi=\Bl1}^{[\bn \cdot \bt]}:=\sum_{i_1=1}^{[n_1t_1]}\dots \sum_{i_d=1}^{[n_dt_d]}.     
\end{equation}
We are looking for the normalization and centering  for $\bbs_{[\bn \cdot \bt]}(\xi)$ in order to get a limit process when $\bn \to \infty$. This means that we consider  the process $\bY(\bt)=\bbs_{[ \bt]}(\xi-\bb)$ (here $\bb$ is some vector of the centering ) and we use the transformation of "time" parameter $\bt$ by diagonal matrix ${\rm diag}(n_1, \dots, n_d)$. Then it is natural for normalization to take diagonal matrices as $F(g)$ , too.

We start by taking $G$ diagonal matrices with positive entries, and since $G$ is homomorphic to $\bbr^d_+$, we shall use the same notation as in Definition \ref{multiselfsim}. Also, as in Definition \ref{multiselfsim}, in Definition \ref{domattr}  and in Theorem \ref{prop4} $T$ stands for $\bbr^d$ or $\bbr^d_+$.

\begin{definition} \label{domattr} A  random field $\{\bY(\bt)=(Y_1(\bt), \dots , Y_m(\bt)), \ \bt \in T\}$, with values in $\bbr^m$, $m\ge 1$, belongs to the domain of  attraction, defined by diagonal matrices with positive entries, of a $\bbr^m$-valued random field $\{\bbx(\bt)=(X_1(\bt), \dots , X_m(\bt)), \ \bt \in T\}$,   (in notation $\bY \in {\rm DA}_{\rm diag}(\bbx)$) , if there exists a    function $\bbf: \bbr_+^d \to \bbr_+^m$  such that, as $\ba\to \infty$,
 \begin{equation}\label{domattreq}
\left \{\frac{\bY(\ba\cdot \bt)}{\bbf(\ba)}, \ \bt \in T\right \}\fdd \{X(\bt), \ \bt \in T \}.
\end{equation}
\end{definition}
 We recall that division of vectors is understood coordinate-wise, also it is possible to consider $\bbf (\ba)$ as diagonal matrix, $\bbf (\ba)= {\rm diag}(f_1(\ba), \dots,f_m(\ba))$, then this ratio would be written as $(\bbf (\ba))^{-1}\bY(\ba\cdot \bt)$.

\begin{teo}\label{prop4} Suppose that a random field $\{\bbx (\bt), \bt \in T\}$ is continuous in probability and $\bY \in {\rm DA}_{\rm diag}(\bbx)$. Then there exist
$${\bbh^{(j)}}=(H_1^{(j)}, \dots, H_d^{(j)}), \ H_i^{(j)}>0,\ 1\le i\le d,$$
 such that
 $$\bbf(\ba)= (f_1(\ba), \dots,f_m(\ba)), \ \ f_j(\ba)=\prod_{i=1}^d a_i^{H_i^{(j)}}L_j(\ba), \ j=1, \dots , m,$$
  where $L_j, \ j=1, \dots , m, $ are  s.v.f.'s, satisfying (\ref{slowvar}), and $\bbx$ is  ${\bar \bbh}$-m.s.s. random field.
\end{teo}

{\it Proof }. The proof of this theorem follows the proof of Theorem \ref{thm2} from \cite{Lamperti}, cited in the Introduction.
 It is clear that the condition (\ref{domattreq}) implies the same condition for each coordinate,  that is, for any $1\le j\le m$,
 \begin{equation}\label{domattreqcor}
\left \{\frac{Y_j(\ba\cdot \bt)}{f_j(\ba)}, \ \bt \in T\right \}\fdd \{X_j(\bt), \ \bt \in T \}.
\end{equation}
We take one fixed  $j, \ 1\le j\le m$, and
 let $\{\ba_i=(a_{i,1}, \dots , a_{i,d}),  \ i\ge 1,\}$ be some sequence satisfying $\ba_i \to \infty$ (we recall that this means $\min_{1\le j\le d} a_{i,j} \to \infty$ as $i\to \infty$). Taking one-dimensional distribution at point $\bt=\Bl1:=(1, \dots , 1)$ from (\ref{domattreqcor}) we get that for almost all $x$
\begin{equation}\label{domattr1}
\PP \{Y_j(\ba_i)\le f_j(\ba_i)x\}\to \PP \{X_j(\Bl1)\le x\}.
\end{equation}
In a similar way, writing $\ba_i=\ba_i \cdot \bt\cdot \bt^{-1}, $ from (\ref{domattreqcor}) we have
\begin{equation}\label{domattr2}
\PP \{Y_j(\ba_i)\le f_j(\ba_i \cdot \bt)x\}\to \PP \{X_j(\bt^{-1})\le x\}.
\end{equation}
Taking into account Lemma in \cite{Lamperti}, p 71, which is a slight extension of the theorem from classical monograph \cite{Gnedenko}, from (\ref{domattr1}) and  (\ref{domattr2}) we get that there exists a limit
$$
\lim_{i\to \infty}\frac{f_j(\ba_i \cdot \bt)}{f_j(\ba_i)}.
$$
Since this limit exists for any  sequence $\{\ba_i\}$,  we have
$$
\lim_{\ba \to \infty}\frac{f_j(\ba \cdot \bt)}{f_j(\ba)}= l_j(\bt).
$$
But this means that $f_j$ is a r.v.f., therefore, applying Proposition \ref{prop2}, we get that $f_j(\bt)=\prod_{i=1}^d t_i^{H_i^{(j)}}L_j(\bt)$ and $l_j(\bt)=\prod_{i=1}^d t_i^{H_i^{(j)}}.$ Thus, we proved the claimed expression for $\bbf$.

It remains to verify that the limit random field is  ${\bar \bbh}$-m.s.s. random field, and this is done exactly as in \cite{Lamperti}. Let us take fixed $k$ points $\bt_i, \ i=1, \dots , k$, then due to (\ref{domattreq})
$$
\left \{\frac{\bY(\ba \cdot \bb \cdot \bt_i )}{\bbf(\bb)}, \ i=1, \dots ,k \right \} \fdd \left \{(\bbx(\ba \cdot \bt_i ), \ i=1, \dots ,k. \right \}
$$
On the other hand, taking into account that
$$
\lim_{\bb\to \infty}\frac{\bbf(\ba \cdot \bb )}{\bbf( \bb)}=\bll(\ba):=\left (l_1(\ba), \dots, l_m(\ba) \right ),
$$
we have
$$
\left \{\frac{\bY(\ba \cdot \bb \cdot \bt_i )}{\bbf(\ba \cdot \bb )}\frac{\bbf(\ba \cdot \bb )}{\bbf(\bb)}, i=1, \dots ,m \right \} \fdd \left \{\bll (\ba)\bbx(\bt_i), i=1, \dots ,m. \right \}
$$
Since $k$ and points $\bt_i, \ i=1, \dots ,k$ can be taken arbitrary, we conclude
$$
 \{\bbx(\ba  \cdot \bt ), \ \bt \in T\}\stackrel{\rm d}{=} \{\bll (\ba)\bbx(\bt), \ \bt \in T\}.
$$
\halmos

This theorem (with $T=\bbr^d_+$)  can be applied to the random field  $\bY(\bt)=\bbs_{[ \bt]}(\xi -b)$, where $\bbs_{[ \bt]}(\xi)$ is  defined in (\ref{sumfields}).

\begin{cor} \label{cor1} Suppose that $\{\xi_{\bi}, \bi\in \bz^d\}$ is a $\bbr^m$-valued stationary random field.
If there exists a function $\bbf(\bn)=(f_1(\bn), \dots,f_m(\bn))\to \infty$ as $\bn\to \infty$ and $\bb\in \bbr^m$ such that
 \begin{equation}\label{coreq}
\left \{\frac{\bbs_{[\bn \cdot \bt]}(\xi -\bb)}{\bbf(\bn)}, \ \bt \in \bbr_+^d\right \}\fdd \{\bbx(\bt), \ \bt \in \bbr_+^d \},
\end{equation}
  where $\bbx$ is a non-degenerate continuous in probability $\bbr^m$-valued random field, then $f_j(\ba)=\prod_{i=1}^d a_i^{H_i^{(j)}}L_j(\ba)$ with some $H_i^{(j)}>0, 1\le i\le d,\ 1\le j\le m,$ and some c.s.v.f.'s $L_j,$ and $\bbx$ is ${\bar \bbh}$-m.s.s. random field.
\end{cor}

In the case $m=1$ we get a result on summation of real valued stationary random field, which can be considered as generalization of Lamperti theorem for stationary sequences, as it is formulated in \cite{Giraitis}.

 Similarly to Definition \ref{domattr}, we can define domain of attraction of a random field using the group of transformations of "time" parameter $\bt$, used in Definition \ref{selfsimop}. As it was explained before this definition, now we can consider only the case $T=\bbr^d$. We recall that we do not use centering.
\begin{definition} \label{domattr3} A $\bbr^m$-valued random field $\{\bY (\bt), \ \bt \in\bbr^d$, belongs to the domain attraction with operator time-scaling  of a  random field $\{\bbx (\bt), \ \bt \in \bbr^d\}$  (in notation $\bY \in {\rm DA}_{op}(\bbx)$), if there exists a $d\times d$ matrix $E\in Q(\bbr^d)$  and   a function $f:r\to {\rm GL}(\bbr^m))$
  such that, as $r\to \infty$,
 \begin{equation}\label{domattreq3}
\left \{f(r)\bY(r^E \bt), \ \bt \in \bbr^d \right \}\fdd \{\bbx(\bt), \ \bt \in \bbr^d \}.
\end{equation}
\end{definition}

These two Definitions \ref{domattr} and \ref{domattr3} differ by the scaling of "time" parameter $\bt$: in Definition \ref{domattr} the group of diagonal matrices with positive entries, homomorphic to $\bbr_+^d,$ is used, while in Definition \ref{domattr3} one parameter group $\{r^E, \ r>0\}$ with a fixed matrix $E$ is applied. Correspondingly, as normalization functions diagonal matrices and more general matrices are used. Another difference is that in Definition \ref{domattr} both cases of $\bbr_+^d$ and $\bbr^d$ are possible, while in Definition \ref{domattr3}  we can consider only the case of $\bbr^d$. We shall see that this difference is important considering Lamperti transformation.

The analog of Lamperti theorem \ref{thm2} in this case is formulated as follows.

\begin{teo}\label{thm4} Suppose that a $\bbr^m$-valued random field $\bbx (\bt), \ \bt \in \bbr^d,$ is proper and continuous in probability, and $\bY \in {\rm DA}_{op}(\bbx)$. Then the random field $\bbx$ is o.s.s..
\end{teo}

{\it Proof}. Since the proof goes along the same lines as the proof of Theorem 5 in \cite{Hudson} in the case of processes, we introduce the same notation. For any integer $k\ge 1$, we  denote $\bT=(\bt_1, \dots , \bt_k), \ \bt_i\in \bbr^d_+, \ ~{1\le i\le k},$ and $ \bbx(\bT)=\left ( \bbx(\bt_1), \dots , \bbx (\bt_k) \right )$. For  linear operators $A$ on $\bbr^m$  and $B$ on $\bbr^d$  $A\bbx(\bT)$ and $B\bT$ stands for $\left ( A\bbx(\bt_1), \dots , A\bbx (\bt_k) \right )$ and $(B\bt_1, \dots , B\bt_k),$ respectively. For any positive integer $k$ and any $r>0$, let $\H_r^k$ stand for the set of all linear operators $H$ on $\bbr^m$ ($m\times m$ matrices ) such that for the matrix $E$ from (\ref{domattreq3}) and all $\bT\in \left ( \bbr^d_+\right )^k$
\begin{equation}\label{domattreq5}
\left \{ \bbx(r^E \bT))\right \}\stackrel{\rm d}{=}\left \{H \bbx(\bT) \right \}.
\end{equation}
The proof that $\bbx$ is o.s.s.  is divided into five steps (in \cite{Hudson} they are formulated as lemmas).

Step 1. It is proved that for all $k$ and  $r>0$  $\H_r^k$ is not empty.

Step 2. It is proved that for all $k$ and  $r>0$  $\H_r^k$ is closed in the set of all linear operators on $\bbr^m$ under the topology induced by the operator norm.

Step 3. It is proved that for all $k$ and  $r>0$  $\H_r^k \supset \H_r^{k+1}.$

Step 4. For all  $r>0$,  $\H_r^1$ is compact in the set of all linear operators on $\bbr^m$.

Step 5. For all  $r>0$,  $\cap_{k=1}^\infty \H_r^k$ is not empty.

The main step, which uses the assumption (\ref{domattreq3}) in Definition \ref{domattr3} and which is connected with centering, is Step 1. Using the associativity of multiplication of matrices we have $(sr)^E\bt=s^Er^E\bt=s^E(r^E\bt)$, therefore, assuming $r>0$,  $k$ fixed, and $s\to \infty$, from (\ref{domattreq3}) we have
\begin{equation}\label{domattreq6}
\left \{f(sr)Y((sr)^E \bT)  \right \}\fdd \{X(\bT) \}
\end{equation}
and
\begin{equation}\label{domattreq7}
\left \{f(s)Y(s^E(r^E \bT)) \right \}\fdd \{X(r^E \bT)\}.
\end{equation}
Since $\bbx$ is proper, it follows that, for sufficiently large $s$, $f(s)$ and $f(sr)$ are invertible and, as in \cite{Hudson}, we can apply Theorem 2.3 from \cite{Weissman} and get that a sequence of linear operators (considered as operators on $(\bbr^m)^k$) $\{ f(n)(f(nr))^{-1}, \ n\ge 1\}$ is relatively compact. If $H$ is a limit point of this sequence, from (\ref{domattreq6}) and (\ref{domattreq7}) we get (\ref{domattreq5}). Thus, $H \in \H_r^k$,  and step 1 is proved. As a matter of fact, the above mentioned Theorem 2.3 from \cite{Weissman} concerns more general affine transformations, i.e., transformations of type $\a:\bbr^m \to \bbr^m, \ \a (\bx)=A\bx+\bh$, where $A$ is a linear operator and $\bh\in \bbr^m$ is fixed. Therefore, in \cite{Hudson} it was possible to use centering, but independent of $t$, see (\ref{domattreq4}). It is easy to see that in our Definition \ref{domattr3}, in the relation (\ref{domattreq3}), we can add centering function $\ba(r)$ and then to prove that $\bbx$ is HMo.s.s.. But if we add centering depending on $\bt$, i.e., instead of (\ref{domattreq3}) we assume the relation
\begin{equation}\label{domattreq8}
\left \{f(r)Y(r^D \bt)+ \ba(r, \bt), \ \bt \in \bbr^d \right \}\fdd \{X(\bt), \ \bt \in \bbr^d \}
\end{equation}
with centering function $\ba: [0, \infty)\times\bbr^d \to \bbr^m$, then we face the following problem. First, we must redefine the set  $\H_r^k$, instead of (\ref{domattreq5}) we must require
$$
\left \{ \bbx(r^D \bT))\right \}\stackrel{\rm d}{=}\left \{H \bbx(\bT)+{\bar \bh} \right \},
$$
 where ${\bar \bh}=(\bh_1, \dots , \bh_k)$ is some collection of vectors from $\bbr^m$.
Second, with centering the relations (\ref{domattreq6}) and (\ref{domattreq7}) become
\begin{equation}\label{domattreq10}
\left \{f(sr)Y((sr)^E \bT)+\ba (sr, \bT)  \right \}\fdd \{X(\bT) \}
\end{equation}
and
\begin{equation}\label{domattreq11}
\left \{f(s)Y(s^E(r^E \bT))+\ba (s, r^E \bT) \right \}\fdd \{X(r^E \bT)\},
\end{equation}
 respectively, where $\ba (s, \bT)=(\ba (s, \bt_1), \dots , \ba (s, \bt_k) )$ for $\bT=(\bt_1, \dots , \bt_k)$. The affine transformation from relation  (\ref{domattreq10}) (taking $\bT$ with one element $\bt$) would be $f(sr)\bx+\ba (sr, \bt)$, therefore, considering $\bt$ as fixed and denoting $\bb(s)=\ba (s, \bt)$ one can take the affine transformation $\a (x)=f(s)x+\bb(s)$. Unfortunately, this transformation is not consistent with the relation (\ref{domattreq11}), since in this relation the affine transformation is $f(s)x +\ba (s, r^E \bt)$ and $\ba (s, r^E \bt)\ne \bb(s)$. At present  we do not know how to overcome these difficulties, arising in the first step of the proof. But we think that there is even no need to do this, since the formulation of the domain of attraction, as it is done in
 (\ref{domattreq4}) and (\ref{domattreq8}) is not natural, since at first there is normalization and then centering, while natural way is opposite - first we use centering and then normalization, as it is, for example, in Theorem \ref{thm2}. Thus, having a process $\bY(\bt)$ we find a centering function $\ba(\bt)$ (if there is need for centering; from examples in summation theory we know that there are situations where centering is not needed) and apply normalization for centered process: $f(r)(\bY(r^E\bt)-\ba(r^E\bt))$. Denoting by ${\bar \bY}(\bt)=\bY(\bt)-\ba(\bt)$ we are in the situation of Definition \ref{domattr3} and Theorem \ref{prop4}.

The proof of the rest four steps can be carried with small changes in the corresponding proofs in \cite{Hudson}.
 Having carried all five steps, it is easy to complete the proof. For any fixed $r>0$, we can take $B(r)\in \cap_{k=1}^\infty \H_r^k$, and, due to (\ref{domattreq5}), we get (\ref{defgroup4}) with $\ba_r(\bt)\equiv 0$.
\halmos

In the case $d=1, E=1$ similar result is proved in \cite{Hudson}. Namely, in Theorem 5 in \cite{Hudson}  instead of (\ref{domattreq3}) the following relation is assumed
\begin{equation}\label{domattreq4}
\left \{f(r)\bY(rt)+ \ba(r), \ t\ge 0 \right \}\fdd \{\bbx(t), \ t\ge 0  \},
\end{equation}
 where $\ba: (0, \infty) \to \bbr^m$ and  $f:r\to {\rm GL}(\bbr^m))$. Then it is  proved that $\bbx$ is HMo.s.s.. Here it is worth to note, that looking at the proof of this result, one can see that in (\ref{domattreq4}) it is possible to consider also the case $t\in R$ (instead of $\{t\ge 0\}$).  It would be more natural to take centering function $\ba$ dependent on $t$, i.e., to take $\ba: (0, \infty)\times\bbr_+ \to \bbr^m$, but then there come these  difficulties which were explained in the proof of  Theorem \ref{thm4}.

\section{Relation between scale transition and domains of attraction}

Recently in the papers \cite{Puplinskaite2} and \cite{Puplinskaite1}    the so-called phenomenon of the scale transition for random fields in the case $d=2$ was described. The case $d>2$ is more complicated and at present the phenomenon of the scale transition in this case remains not investigated.  It is interesting how  this phenomenon is related to the analog of Lamperti theorem, formulated in Corollary \ref{cor1} with $m=1$. The scaling transition phenomenon can be explained  briefly as follows. For real-valued random fields in the case $d=2$ we shall use different notation. Suppose that we consider a stationary random field $\xi=\{\xi_{i,j}, \ (i,j)\in \bz^2\}$ and sum-processes
\begin{equation}\label{scaletr}
S_{n, m}(t,s)=\sum_{i=1}^{[nt]} \sum_{j=1}^{[m s]} \xi_{i,j}, \quad t\ge 0, \quad s\ge 0, \quad Z_{n, \g}(t,s)=S_{n, n^{\g}}(t,s)
\end{equation}
where $\g >0$ is a parameter, and $Z_{n, \g}(t,s)=0$, if at least one of $[nt]$ or $[n^\g s]$ is zero. It is assumed that, for any $\g>0$, there exists a nontrivial random field $V_\g (t,s)$ and a normalization $A_n (\g) \to \infty$ such that f.d.d. of $A_n^{-1} (\g)Z_{n, \g}(t,s)$ converges weakly to f.d.d. of $V_\g$. It is said that the random field $\xi$ exhibits scaling transition if there exists $\g_0>0$ such that the limit process $V_\g$ is the same, let say $V_+$, for all  $\g >\g_0$ and  another, not obtained by simple scaling, $V_-$, for   $\g <\g_0$. In  \cite{Puplinskaite2} it is demonstrated that  for a long-range dependent Gaussian random field $\xi$ under some conditions  on spectral density of the field this phenomenon of the scale transition can be observed. In this case we face a problem when considering the convergence of f.d.d. of random fields $S_{n, m},$ defined in (\ref{scaletr}), since the limit process depends on the way how  $(n,m)$ tends to infinity: if we take $m=n^\g$ then we have that $\lim_{n \to \infty}m/n$ can be  $0$ (if $\g <1$), $1$ (if $\g=1$), or $\infty$ (if $\g >1$), and the limit process can be different for different values of $\g$.

Condition (\ref{coreq}) from Corollary \ref{cor1} in our notation can be rewritten as follows: there exists a function $f(n,m)\to \infty$, as $(n, m)\to \infty$, such that
\begin{equation}\label{coreq1}
\left \{\frac{S_{n, m}(t,s)}{f(n, m)}, \ (t, s) \in \bbr_+^2\right \}\fdd \{X(t, s), \ (t, s) \in \bbr_+^2 \},
\end{equation}
 where $X$ is a non-degenerate continuous in probability real-valued random field. Thus, there is a
 requirement that the limit process in (\ref{coreq1}) is independent on the way how $(n,m)$ tends to infinity. It means that Corollary \ref{cor1} cannot be applied to those random fields exhibiting scale transition. On the other hand, we may restrict the ways  how  $(n,m)$ tends to infinity, namely,
 to the requirement $\min(m,n)\to \infty$ we may add the following condition: there exist constants $0<c<C<\infty$ such that
\begin{equation}\label{limit1}
c \le \liminf \frac{m}{n}\le \limsup \frac{m}{n}\le C.
\end{equation}
Such condition, which is quite natural, will exclude the possibility of the scale transition. In  case  $d>2$ in order to exclude scale transition (although at present the case $d>2$ remains not investigated, but it is clear  that this phenomenon will be, and in higher dimensions it will be  even more complicated) we must control all ratios $n_i/n_j, \ i=1, \dots, d-1, \ j>i,$ \ in a similar way as in (\ref{limit1}). That  condition (\ref{limit1}) is natural shows the following argument. If we  additionally suppose that there exits limit $\lim_{n\to \infty} m/n =c$, then scaling transition phenomenon is possible only if $c=0$ or $c=\infty$, while for $0<c<\infty$ condition (\ref{limit1}) holds. One of areas, where two-dimensional data sets arise, is analysis of panel data
$\{\xi_{i,j}, \ 1\le i\le n, \ 1\le j \le T,\}$, where $n$ is the total number of panels (or cites where data is collected from) and T is total number of data collected at moments $t_j, 1\le j  \le T$ at each panel. Various limit theorems (under different assumptions on data) for sums $\sum_{i=1}^{n} \sum_{j=1}^{T} \xi_{i,j}$  are considered in many papers, we shall point here only the paper \cite{PhillMoon}, where two types of convergence of the above written sums are considered: the sequential convergence (at first $T\to \infty$, then $n\to \infty$, and the order of these limits is important),
joint convergence ($(n, T) \to \infty$ jointly) and the so called diagonal convergence when $T=T(n)$ (the last type of convergence is not considered in the cited paper). Let us note that in \cite{PhillMoon} the joint convergence means that $n\to \infty, \ T\to \infty$, but also condition (\ref{limit1}) holds, although this condition is not stated explicitly. This can be understood since in several theorems it is written: "$(n, T)\to \infty$ with $n/T \to 0$". If $(n, T)\to \infty$ would mean, as in our paper only $\min (n, T)\to \infty$ (which clearly includes cases $0\le \lim n/T \le \infty$), there would be no sense to add  $n/T \to 0$. The only diagonal convergence may lead to the scale transition phenomenon, if we assume $T=n^\g$ with $\g>1$, but it seems that in panel data analysis, as it is written in \cite{PhillMoon}, the diagonal case is understood only as  $n/T \to c\ne 0$. Thus, for practitioner, having data with $n=82$ and $T=10000$, there are two possibilities: to apply the sophisticated theorem, obtained assuming $T=n^\g$ with $\g>1$, but not knowing if there is scale transition and if it is, what is the value of the critical point $\g_0$, or to apply theorem obtained under assumption $n\le cT$ with small $c$ (let us say, $c=10^{-3}$). Most probably, a practitioner will choose the latter possibility.

\section{Lamperti type transformation}

Finally, we discuss the generalization of Proposition \ref{prop1} for self-similar random fields. Having general Definitions \ref{groupss}  and \ref{rdgroupss} of self-similarity, the Lamperti transformation, given in Proposition \ref{prop1}, can be considered as the relation between two types of self-similarity, since stationarity can be considered as self-similarity with respect to the group of translations. Therefore, trying to find Lamperti transformation for general $(G, C)$-ss random fields from Definition \ref{groupss}, only with cocycle $C$ independent of $t$, at first we must formulate abstract analogue of stationarity. Let $S$ be some set and let $\H$ be a group of transformations of $S$.
\begin{definition} \label{defsimilar} We say that a process $\{ \bY(s),  s\in S\}$ is $\H$-stationary if for all $h\in \H$
\begin{equation}\label{similar}
\left \{ \bY(h(s)), \ s\in S)\right \}\lygfdd \left \{ \bY(s), \ s\in S)\right \}.
\end{equation}
\end{definition}
 Now we want to establish the relation between $(G, C)$-ss random field $\bbx$ and $\H$-stationary random field $\bY$. Clearly, we must assume that groups of transformations $G$ and $\H$, acting on $T$ and $S$, respectively, must be compatible in some sense, and we assume that $G$ and $\H$ are isomorphic and $F: \H \to G$ (bijective function preserving corresponding operations in groups $G$ and $\H$) represents this isomorphism. Also let us denote by $\C$ the group of bijections in $\bbr^m$, generated by cocycle $C$ and a group  $G$, i.e., $\C=\{C(g), g\in G\}$. Now the problem can be formulated as follows. Let $\{\bbx(t), \ t\in T\}$ satisfies the relation
 \begin{equation}\label{defgroup1a}
\left \{ \bbx(g(t)), \ t\in T \right \}\stackrel{\rm d}{=}\left \{C(g) \bbx(t), \ t\in T \right \},
\end{equation}
i.e., $\bbx$ is $(G, C)$-ss. We would like to find a set $S$,  a group of  transformations $\H$ of $S$, an invertible mapping $\vfi:S\to T$, and a mapping $f:S\to \C$ such, that the process
\begin{equation}\label{stacionar}
\bY(s)= f(s) \bbx(\vfi (s)), \ s\in S,
\end{equation}
would be $\H$-stationary. We assume that the mappings $F$ and $\vfi$ commute in the following sense: for all $s\in S, \ h\in \H$,
 \begin{equation}\label{cond2}
F(h)(\vfi (s))=\vfi (h(s)).
\end{equation}

 Now we can formulate the following generalization of Proposition \ref{prop1}.

\begin{prop}\label{prop6}. Let  $\bbx$ be $(G, C)$-ss. If the mappings $\vfi, f, $ and $F$ are as they are introduced above and, for all $s\in S, \ h\in \H$, the following relation
\begin{equation}\label{cond1}
f(h(s))C(F(h))=f(s)
\end{equation}
 holds, then the process $\{ \bY(s),  s\in S\}$, defined in (\ref{stacionar}) is $\H$-stationary.

 Conversely, let $\{ \bY(s),  s\in S\}$ be $\H$-stationary and let $T, G,$ and $C$ be given. If the mappings $F, \vfi, f$ satisfy (\ref{cond2}) and (\ref{cond1}), then the process
 \begin{equation}\label{converse}
\bbx (t)=\left (f(\vfi^{-1}(t))\right )^{-1}\bY (\vfi^{-1}(t))
\end{equation}
is $(G, C)$-ss.
\end{prop}

{\it Proof }. Using (\ref{defgroup1a})-(\ref{cond2}) we have
$$
\bY(h(s))=f(h(s))\bbx(\vfi (h(s)))=f(h(s))\bbx(F(h)(\vfi (s)))\stackrel{\rm d}{=}f(h(s))C(F(h))\bbx(\vfi (s)).
$$
From this relation, using (\ref{cond1}), we get (\ref{similar}).

To prove the converse relation we must show that the process defined in (\ref{converse}) satisfies the relation
\begin{equation}\label{defgroup7}
\left \{ \bbx(g(t)), t\in T \right \}\stackrel{\rm d}{=}\left \{C(g) \bbx(t), t\in T \right \}.
\end{equation}
We have
$$
\bbx(g(t))=\left (f(\vfi^{-1}(g(t)))\right )^{-1}\bY (\vfi^{-1}(g(t))).
$$
Remembering that  $\vfi$ is invertible, it is easy to see that  (\ref{cond2}) implies the following relation: for all $t\in T, \ h\in \H$, if $F(h)=g$, then
\begin{equation}\label{cond3}
\vfi^{-1}g(t)=h(\vfi^{-1}(t)).
\end{equation}
Using stationarity of the process $\bY$ and (\ref{cond3}) we can write
$$
\bY (\vfi^{-1}(g(t)))=\bY (h(\vfi^{-1}(t)))=\bY (\vfi^{-1}(t))=\bY(s)
$$
and
$$
\left (f(\vfi^{-1}(g(t)))\right )^{-1}=\left (f(h(\vfi^{-1}(t))\right )^{-1}=\left (f(h(s))\right )^{-1}f(s)(f(s))^{-1}.
$$
From (\ref{cond1}), remembering that $F(h)=g$, we have
$$
C(g)=\left (f(h(s))\right )^{-1}f(s).
$$
Collecting the obtained relations, we get
$$
\bbx(g(t))=\left (f(h(s))\right )^{-1}f(s)(f(s))^{-1}\bY(s)=C(g)\bbx (t),
$$
and (\ref{defgroup7}) is proved.
\halmos

Before providing simple corollaries from this result, let us clarify the meaning of the condition (\ref{cond1}). Suppose that we have $h_0\in \H, \ s_0\in S$ and $s_1=h_0(s_0), \ F(h_0)=g_0$. Since $f$ maps $S$ to  $\C=\{C(g), g\in G\}$, let us denote $f(s_0)=C(g_1), \ f(s_1)=C(g_2).$ Writing  (\ref{cond1}) for $h_0$ and $s_0$ we have
$$
f(h_0(s_0))C(g_0))=f(s_0),
$$
but this relation gives us the following equality
$$
C(g_2)C(g_0)=C(g_1).
$$
Since $C$ is a cocycle, this means that $F$ and $f$ must be such that $g_1=g_2g_0$.

We can provide one case, where we can verify the condition (\ref{cond1}). Suppose that $S$ itself has structure of a group, then we can identify $\H$
with $S$ and to show that
\begin{equation}\label{fs1}
f(s)=\left (C(F(s))\right )^{-1}
\end{equation}
satisfies (\ref{cond1}). Since now  $F: S\to G$ and $F(s)(F(s))^{-1}$ is identical map on  $G$, let say  $e$, and $C(e)=I_m$ (identity on $\bbr^m$), therefore, applying (\ref{cocycle})  we get
\begin{equation}\label{fs2}
\left (C(F(s))\right )^{-1}=C((F(s))^{-1}).
\end{equation}
Now we can write, taking into account (\ref{fs2}) and still using the notation $h(s)$, despite the fact that now both $h$ and $s$ are elements of $S$,
$$
f(h(s))=\left (C(F(s))\right )^{-1}=C((F(h(s)))^{-1})
$$
$$
=C\left ((F(h)F(s))^{-1}\right )=C\left ((F(s))^{-1}\right )\cdot C\left ((F(h))^{-1}\right ).
$$
Therefore,
$$
f(h(s))C(F(h))=C\left ((F(s))^{-1}\right )\cdot C\left((F(h))^{-1}\right )\cdot C(F(h))=C\left ((F(s))^{-1}\right )=f(s),
$$
that is , we get (\ref{cond1}).

The classical Lamperti transformation is obtained from Proposition \ref{lamperti-rm} by taking $T=G=\bbr_+, \ C(g)(x)=g^{H}x,  \H=S=\bbr, \ h(s)=s+h, \ f(s)=e^{-H s}, \ \vfi (s)=e^{s}, \ F(h)=e^{h}$. Then it is easy to see that $Y(s)=e^{-H s}X(e^s)$ is stationary, if $X$ is $H$-ss. Conversely, if $\{Y(s), \ s\in \bbr\}$ is a stationary process, then $X(t)=t^HY(\ln t), t\in \bbr_+$ is $H$-ss.

Lamperti transformation for $(\bbh)-mss$ real valued random fields was given in  \cite{Genton} (see Proposition 1 therein ),  from Proposition \ref{prop6} we shall get more general result for $\bbr^m$-valued multi-self-similar random fields, defined in Definition \ref{multiselfsim}.

 We recall that  ${\bar \bbh}$-m.s.s.   stands for $\bbr^m$-valued multi-self-similar  random field. From Proposition \ref{prop6} we get the following result.

\begin{prop}\label{lamperti-rm} Let $\{\bbx(\bt), \ \bt \in \bbr^d_+\}$ be ${\bar \bbh}$-m.s.s. random field. Then the $\bbr^m$ process
\begin{equation}\label{stacionar1}
\bY(\bs)=\left \{\exp \left (-\sum_{i=1}^d s_iH_i^j\right )X_j(e^{s_1}, \dots, e^{s_d} ), \ j=1, \dots, m, \bs \in \bbr^d  \right \}
\end{equation}
is stationary with respect to the usual translation operation on $\bbr^d$. Conversely, if $\bY(\bt), \bt \in \bbr^d,$ is stationary, then
\begin{equation}\label{Hmss}
\bbx(\bt)=\left \{\prod_{i=1}^d t_i^{H_i^j}Y_j(\ln {t_1}, \dots, \ln{t_d} ), \ j=1, \dots, m, \ \bt \in \bbr^d_+  \right \}
\end{equation}
is  ${\bar \bbh}$-m.s.s. random field.
\end{prop}

{\it Proof }. To derive this proposition from Proposition \ref{prop6} we take $T= \bbr^d_+, S=\bbr^d$, and we   identify  $\H$  with $S$. Then we take  $\vfi(\bs)=(e^{s_1}, \dots, e^{s_d} ), \ h(\bs)=\bs+\bh, \ \bs,\bh \in \bbr^d$,
$$
F: \bbr^d \to \bbr^d_+, \ F(\bh) ={\rm diag}(e^{h_1}, \dots, e^{h_d} ), \quad (F(\bh))^{-1} ={\rm diag}(e^{-h_1}, \dots, e^{-h_d} ),
$$
$$
F^{-1}: \bbr^d_+ \to \bbr^d, \ F^{-1}(\bs) =(\ln s_1, \dots , \ln s_d), \ \bs \in \bbr^d_+
$$
$$
C(g)={\rm diag}\left (\prod_{i=1}^d a_i^{H_i^j}, \ j=1, \dots, m \right ), \quad {\rm for} \quad g={\rm diag}(a_1, \dots , a_d).
$$
It is easy to verify that the above defined functions satisfy the relations (\ref{cond2}) and (\ref{cond1}).
Since we    identified $\H$  with $S$, we can use (\ref{fs1}) and (\ref{fs2}) and get
$$
f(\bs)=C((F(\bs))^{-1})={\rm diag}\left (\prod_{i=1}^d s_i^{H_i^j}, \ j=1, \dots, m \right ).
$$
Now simple substitution  of the obtained expressions into (\ref{stacionar}) and (\ref{converse}) allows to get (\ref{stacionar1}) and (\ref{Hmss}).
\halmos

More complicated situation is with o.s.s. random fields. In \cite{Hudson}  the following result for processes  was formulated.
\begin{prop}\label{lamperti-HMrm} (\cite{Hudson})
 Suppose that  $\bbr^m$-valued random process $\{\bbx (t),  t\ge 0\}$ is proper and  satisfies the relation
\begin{equation}\label{ossLamp}
\left \{ \bbx(rt), t\ge 0\right \} \stackrel{\rm d}{=}\left \{r^D \bbx(t), t\ge 0\right \},
\end{equation}
then $\{\bY (t):=e^{-tD}\bbx(e^t), \ t\in \bbr\}$ is a stationary process. Conversely, let $\{\bY(t), \ t\in \bbr\}$ be a stationary process and let $D$ be a nonsingular $m\times m$ matrix such that $t^D\bY(\ln t)$, as $t\to 0$, converges in law to some distribution $\mu$. Then  the process $\{\bbx (t), \ t\ge 0\},$ defined by $\bbx (t):=t^D\bY(\ln t)$, for $t>0$, and for $t=0$ as  a random variable with distribution $\mu$ is proper and o.s.s., satisfying (\ref{ossLamp}).
\end{prop}

It turns out that the  generalization of this result for o.s.s.  random fields is not an easy task. In all cases, where we were able to find Lamperti type transformation, we had $T=\bbr_+$ or $\bbr_+^d$ and
$S=\bbr$ or $\bbr^d$, respectively, and exponential map realized homomorphism between these spaces. Now, dealing with o.s.s.  (or w.o.s.s. , or HMo.s.s.) random fields, we have $T=\bbr^d$, and it is not clear what space $S$ and what bijection $\vfi:S\to T$ can be taken. At present, as far as we know, there is only one result in the case $ m=1$ and $T=\bbr^2$ for particular case of L{\'{e}}vy fractional Brownian random field. Since in this result $S$ is taken also $\bbr^2\setminus \{\B0\}$, but with polar coordinates, for $\bt \in \bbr^2, \ \bt\ne \B0,$ let us denote $\rho(\bt)={\sqrt {t_1^2+t_2^2}}, \ \theta (\bt)={\arctan (t_2/t_1)}+k\pi, k\in \bz.$

\begin{teo}\label{lamperti-levy} (\cite{Genton}) Let $\{X(\bt), \ \bt \in \bbr^2\}$ be a mean zero L{\'{e}}vy fractional Brownian random field with covariance
$$
EX(\bt)X(\bu)=\frac{1}{2}\left (||\bt||^{2H}+||\bu||^{2H}+||\bt-\bu||^{2H} \right ),
$$
where $0<H\le 1$. Then $X(\bt)\lygfdd \rho(\bt)^H Y(\ln \rho(\bt), \theta (\bt))$, where $\{Y(\bs), \ \bs \in \bbr^2\}$ is a mean zero Gaussian stationary process with covariance
\begin{equation}\label{covar}
R(\bv):=EY(\bs)Y(\bs +\bv)=\frac{1}{2}\left (e^{v_1H}+e^{-v_1H}- (e^{v_1}+e^{-v_1}-2\cos (v_2))^H\right ).
\end{equation}
Conversely, if $\{Y(\bt), \bt \in \bbr^2\}$ is a mean zero Gaussian stationary process with covariance $R(\bv)$ given by (\ref{covar}) and $0<H\le 1$, then $Y(\bt)\lygfdd e^{-t_1H}X\left (e^{t_1}\cos t_2,  e^{t_1}\sin t_2\right )$, where $\{X(\bt), \ \bt \in \bbr^2\}$ is a mean zero L{\'{e}}vy fractional Brownian random field.
\end{teo}

But this  is very particular result, even the analog of this result for a L{\'{e}}vy fractional Brownian random field defined on $\bbr^d$ with $d\ge 3$, as it is mentioned in \cite{Genton}, is not known.
More difficult problem is to find Lamperti transformation for ${\bar \bbh}$-m.s.s. random field, defined on $\bbr^d$ (in Proposition \ref{lamperti-rm} the case of $\bbr^d_+$ is considered).

 Considering  an o.s.s.  random field $\bbx (\bt)$ which satisfies the relation
$$
\left \{\bbx(r^E \bt), \bt\in \bbr^d\right \} \stackrel{\rm d}{=}\left \{r^D \bbx(\bt), \bt\in \bbr^d \right \},
$$
one can think that the so-called non-Euclidean polar coordinates, connected with matrix $E$ (see, for example, \cite{Jurek}  or \cite{MeerschaertB}), can be of some help. But even taking $\bbr^d \setminus \{\B0\}$ with these polar coordinates as $S$, one must guess the group of transformations $\H$ on $S$, which, on one hand, would be appropriate for defining stationarity, on the other hand, would be isomorphic with the group $G=\{r^E, \ r>0\}$. This is a difficult problem for future research.

\end{document}